\documentclass[12pt,reqno]{amsart}
\usepackage{amssymb}
\usepackage{amsxtra}
\usepackage{mathrsfs}
\usepackage[all]{xy}
\usepackage{cite}
\usepackage{paralist}
\usepackage[hypertex]{hyperref}
\newtheorem{theorem}{Theorem}[section]
\newtheorem{lemma}[theorem]{Lemma}
\newtheorem{prop}[theorem]{Proposition}
\newtheorem{corollary}[theorem]{Corollary}
\theoremstyle{definition}
\newtheorem{definition}[theorem]{Definition}
\theoremstyle{remark}
\newtheorem{example}[theorem]{Example}
\newtheorem*{examples*}{Examples}
\newtheorem{remark}[theorem]{Remark}
\newtheorem*{ackn}{Acknowledgements}
\DeclareMathOperator{\Hol}{Hol}
\DeclareMathOperator{\Pol}{Pol}
\newcommand*{\Ptens}{\mathop{\widehat\otimes}}
\newcommand*{\Tens}{\mathop{\otimes}}
\newcommand*{\Pfree}{\mathop{\widehat *}}

\newcommand*{\op}{\mathrm{op}}
\newcommand*{\reg}{\mathrm{reg}}
\newcommand*{\hol}{\mathrm{hol}}
\newcommand*{\free}{\mathrm{free}}
\newcommand*{\wh}{\widehat}

\newcommand*{\la}{\langle}
\newcommand*{\ra}{\rangle}
\newcommand*{\CC}{\mathbb C}
\newcommand*{\N}{\mathbb N}
\newcommand*{\Z}{\mathbb Z}
\newcommand*{\R}{\mathbb R}
\newcommand*{\DD}{\mathbb D}
\newcommand*{\BB}{\mathbb B}
\newcommand*{\cO}{\mathscr O}
\newcommand*{\cA}{\mathscr A}
\newcommand*{\cC}{\mathscr C}
\newcommand*{\cF}{\mathscr F}
\newcommand*{\cN}{\mathscr N}
\newcommand*{\cB}{\mathscr B}

\newcommand*{\eps}{\varepsilon}
\newcommand*{\ol}{\overline}
\newenvironment{mycompactenum}{\pltopsep=5pt\begin{compactenum}[\upshape (i)]}%
{\end{compactenum}}
\begin{document}
\title[Noncommutative analogues of Stein spaces]{Noncommutative analogues
of Stein spaces\\ of finite embedding dimension}
\subjclass[2010]{46L89, 46L52, 46H30, 14A22, 16S38, 46L65, 16S80}
\author{A. Yu. Pirkovskii}
\address{Faculty of Mathematics\\
National Research University Higher School of Economics\\
Vavilova 7, 117312 Moscow, Russia}
\email{aupirkovskii@hse.ru, pirkosha@gmail.com}
\date{}

\begin{abstract}
We introduce and study {\em holomorphically finitely generated} (HFG) Fr\'echet algebras,
which are analytic counterparts of affine (i.e., finitely generated)
$\CC$-algebras. Using a theorem of O.~Forster, we prove that
the category of commutative HFG algebras is anti-equivalent to
the category of Stein spaces of finite embedding dimension.
We also show that the class of HFG algebras is stable under some standard
constructions. This enables us to give a series of concrete examples
of HFG algebras, including Arens-Michael envelopes of affine algebras
(such as the algebras of
holomorphic functions on the quantum affine space and on the quantum torus),
the algebras of holomorphic functions
on the free polydisk, on the quantum polydisk, and on the quantum ball.
We further concentrate on the algebras of holomorphic functions
on the quantum polydisk and on the quantum ball and show that they are
isomorphic, in contrast to the classical case. Finally, we interpret
our algebras as Fr\'echet algebra deformations of the classical algebras
of holomorphic functions on the polydisk and on the ball in $\CC^n$.
\end{abstract}

\maketitle

\section{Introduction}

To motivate our constructions, let us start by recalling some basic principles
of noncommutative geometry. For more details, see, e.g., \cite{Khal_bas_NCG}.

Noncommutative geometry is based on the observation that, very often,
all essential information about a geometric space (e.g., an affine
algebraic variety, or a smooth manifold, or a topological space,
or a measure space) is contained in a suitably chosen algebra of
functions on the space. Sometimes this observation leads to
a theorem which establishes an anti-equivalence between a certain
category of spaces and the respective category of algebras of functions on such
spaces. Moreover, the resulting category of algebras often admits an
abstract characterization not involving any functions on any spaces.
So the ideal starting point for any kind of noncommutative geometry
is roughly as follows. Let $\cA$ be a category of associative
algebras (maybe endowed with an additional structure),
and suppose that the full subcategory of $\cA$ consisting
of commutative algebras is anti-equivalent to a certain category
$\cC$ of ``spaces''.
In such a situation we may think of algebras
belonging to $\cA$ as noncommutative analogues of spaces belonging to $\cC$.
Let us give some concrete illustrations of this phenomenon.

\begin{example}
\label{ex:GN}
By the Gelfand-Naimark Theorem, the category of compact Hausdorff
topological spaces is anti-equivalent to the category of commutative unital
$C^*$-algebras via the functor taking a compact topological space $X$
to the algebra $C(X)$ of continuous functions on $X$.
This leads to the idea that arbitrary (i.e., not necessarily
commutative) unital $C^*$-algebras may be viewed as ``noncommutative
compact topological spaces''. The above point of view has proved to
be very productive. It gave birth to noncommutative topology,
a deep and important subject having many applications (see, e.g.,
\cite{Connes_NCG,Connes_Marcolli,Gracia,Khal_bas_NCG}).
\end{example}

\begin{example}
\label{ex:Null}
Hilbert's Nullstellensatz implies that the category of affine algebraic
varieties (over $\CC$) is anti-equivalent to the category of commutative unital
finitely generated algebras without nilpotents. As above, the equivalence
takes an affine variety $V$ to the algebra $\cO^\reg(V)$ of regular functions
on $V$. More generally, the same functor yields an anti-equivalence between
the category of affine schemes of finite type over $\CC$ and the category
of commutative unital finitely generated algebras (with nilpotents allowed).
Thus unital finitely generated algebras may be viewed as ``noncommutative
affine algebraic varieties'', or ``noncommutative affine schemes of finite type''.
This point of view leads to noncommutative affine algebraic geometry, also
a challenging and rapidly growing subject
(see, e.g., \cite{VOV,LeBruyn_NCG,VO_AG_ass,Rosen_NCG,comp_stab_shv}).
\end{example}

\begin{example}
\label{ex:For_thm}
One more illustration is a theorem by O.~Forster. Recall (see, e.g., \cite{GR_II})
that for each complex space $X=(X,\cO_X^\hol)$ the algebra $\cO^\hol(X)$ has a canonical
topology making it into a nuclear Fr\'echet algebra. If $X$ is reduced
(e.g., if $X$ is a complex manifold), then the elements of $\cO^\hol(X)$ are
holomorphic functions on $X$, and the canonical topology on $\cO^\hol(X)$
is the topology of compact convergence. By definition, a Fr\'echet algebra $A$
is a {\em Stein algebra} if it is topologically isomorphic to $\cO^\hol(X)$
for some Stein space $X$. Forster's theorem \cite{For} states that the
functor $(X,\cO^\hol_X)\mapsto\cO^\hol(X)$ is an anti-equivalence between
the category of Stein spaces and the category of Stein algebras.
\end{example}

By comparing the above examples, we see that Example~\ref{ex:For_thm} does not entirely
match the above ``ideal'' setup for noncommutative geometry.
Indeed, in contrast to Examples~\ref{ex:GN} and~\ref{ex:Null}, Forster's theorem
does not give an abstract characterization of Stein algebras, and so it
does not give us any hint of what noncommutative Stein spaces are.

Our goal here is to introduce a category $\cA$ of Fr\'echet algebras such that
the commutative part of $\cA$ is anti-equivalent to the category of Stein spaces
of finite embedding dimension. We hope that such a category may be useful for
developing noncommutative complex analytic geometry, a field
much less investigated than other types of noncommutative geometry (like, for example,
noncommutative algebraic geometry, noncommutative differential geometry,
noncommutative topology, and noncommutative measure theory).

This paper is mostly a survey. The proofs are either sketched or omitted. We plan to present
the details elsewhere.

\section{Holomorphically finitely generated algebras}

We shall work over the field $\CC$ of complex numbers. All algebras are assumed to
be associative and unital, and all algebra homomorphisms are assumed to be unital
(i.e., to preserve identity elements).
By a {\em Fr\'echet algebra} we mean a complete metrizable locally convex
algebra (i.e., a topological algebra whose underlying space
is a Fr\'echet space). A {\em locally $m$-convex algebra} \cite{Michael} is a topological
algebra $A$ whose topology can be defined by a family of submultiplicative
seminorms (i.e., seminorms $\|\cdot\|$ satisfying $\| ab\|\le \| a\| \| b\|$
for all $a,b\in A$). A complete locally $m$-convex algebra is called
an {\em Arens-Michael algebra} \cite{X2}.

Recall that for each Arens-Michael algebra $A$ and for each commuting $n$-tuple
$a=(a_1,\ldots ,a_n)\in A^n$ there exists an {\em entire functional calculus},
i.e., a unique continuous homomorphism $\gamma_a$
from the algebra $\cO(\CC^n)$ of holomorphic functions on $\CC^n$ to $A$
taking the complex coordinates $z_1,\ldots ,z_n$ to $a_1,\ldots ,a_n$, respectively.
Explicitly, $\gamma_a$ is given by
\[
\gamma_a(f)=f(a)=\sum_{\alpha\in\Z_+^n} c_\alpha a^\alpha
\quad\text{for } f=\sum_{\alpha\in\Z_+^n} c_\alpha z^\alpha\in\cO(\CC^n),
\]
where we use the standard notation $a^\alpha=a_1^{\alpha_1}\cdots a_n^{\alpha_n}$
for $\alpha=(\alpha_1,\ldots ,\alpha_n)\in \Z_+^n$. Below we will need a
noncommutative (or, more exactly, free) version of the entire functional
calculus due to J.~L.~Taylor \cite{T2}.

Let $F_n=\CC\la \zeta_1,\ldots ,\zeta_n\ra$ be the free algebra with generators
$\zeta_1,\ldots ,\zeta_n$. For each $k\in \Z_+$, let $W_{n,k}=\{ 1,\ldots ,n\}^k$,
and let $W_n=\bigsqcup_{k\in\Z_+} W_{n,k}$. Thus a typical element of $W_n$
is a $k$-tuple $\alpha=(\alpha_1,\ldots ,\alpha_k)$ of arbitrary length
$k\in\Z_+$, where $\alpha_j\in\{ 1,\ldots ,n\}$ for all $j$. The only element
of $W_{n,0}$ will be denoted by $\emptyset$. For each
$\alpha=(\alpha_1,\ldots ,\alpha_k)\in W_n$ we let
$\zeta_\alpha=\zeta_{\alpha_1}\cdots\zeta_{\alpha_k}\in F_n$ if $k>0$;
it is also convenient to set $\zeta_\emptyset=1\in F_n$. The set
$\{\zeta_\alpha : \alpha\in W_n\}$ of all words in $\zeta_1,\ldots ,\zeta_n$
is the standard vector space basis of $F_n$.

For each $\alpha\in W_{n,k}\subset W_n$, let $|\alpha|=k$.
The algebra of {\em free entire functions} \cite{T2,T3} is defined to be
\begin{equation}
\label{Fn}
\cF_n=\Bigl\{ a=\sum_{\alpha\in W_n} c_\alpha\zeta_\alpha :
\| a\|_\rho=\sum_{\alpha\in W_n} |c_\alpha|\rho^{|\alpha|}<\infty
\;\forall\rho>0\Bigr\}.
\end{equation}
The topology on $\cF_n$ is given by the seminorms $\|\cdot\|_\rho\; (\rho>0)$,
and the product is given by concatenation (like on $F_n$).
Each seminorm $\|\cdot\|_\rho$ is easily seen to be submultiplicative, so
$\cF_n$ is a Fr\'echet-Arens-Michael algebra containing $F_n$ as a dense
subalgebra. As was observed by D.~Luminet \cite{Lum}, $\cF_n$ is nuclear.
Note also that $\cF_1$ is topologically isomorphic to $\cO(\CC)$.

Taylor \cite{T2} observed that for each Arens-Michael algebra $A$
and for each $n$-tuple
$a=(a_1,\ldots ,a_n)\in A^n$ there exists a {\em free entire functional calculus},
i.e., a unique continuous homomorphism $\gamma_a^\free\colon\cF_n\to A$
taking the free generators $\zeta_1,\ldots ,\zeta_n$ to $a_1,\ldots ,a_n$, respectively.
Explicitly, $\gamma_a^\free$ is given by
\begin{equation}
\label{free_entire}
\gamma_a^\free(f)=f(a)=\sum_{\alpha\in W_n} c_\alpha a_\alpha
\quad\text{for } f=\sum_{\alpha\in W_n} c_\alpha \zeta_\alpha\in\cF_n,
\end{equation}
where $a_\alpha=a_{\alpha_1}\cdots a_{\alpha_k}\in A$ for
$\alpha=(\alpha_1,\ldots ,\alpha_k)\in W_n$, and $a_\emptyset=1\in A$.

From now on, we assume that $A$ is a Fr\'echet-Arens-Michael algebra.

\begin{definition}
We say that a subalgebra $B\subseteq A$ is {\em holomorphically closed}
if for each $n\in\N$, each $b\in B^n$, and each $f\in\cF_n$ we have
$f(b)\in B$.
\end{definition}

In the commutative case, we may use $\cO(\CC^n)$ instead of $\cF_n$:

\begin{prop}
Let $A$ be a commutative Fr\'echet-Arens-Michael algebra.
Then a subalgebra $B\subseteq A$ is holomorphically closed if and only
if for each $n\in\N$, each $b\in B^n$, and each $f\in\cO(\CC^n)$ we have
$f(b)\in B$.
\end{prop}

\begin{examples*}
Of course, if $B$ is closed in $A$, then it is holomorphically closed.
More generally, if $B$ is an Arens-Michael algebra under a topology
stronger than the topology inherited from $A$, then $B$ is holomorphically
closed in $A$. For example, the algebra $C^\infty(M)$ of smooth functions
on a compact manifold $M$ is holomorphically closed in $C(M)$.
\end{examples*}

\begin{remark}
Our notion of a holomorphically closed subalgebra should not be confused with
the more common notion of a subalgebra stable under the holomorphic functional
calculus (see, e.g., \cite{Connes_NCG,Gracia}). Recall that a subalgebra
$B\subseteq A$ is stable under the holomorphic functional
calculus if for each $b\in B$, each neighbourhood $U$ of the spectrum $\sigma_A(b)$,
and each $f\in\cO(U)$ we have $f(b)\in B$. Such a subalgebra is
necessarily spectrally invariant in $A$, i.e., $\sigma_B(b)=\sigma_A(b)$ for all
$b\in B$. In contrast, a holomorphically closed subalgebra need not be spectrally
invariant. For example, for each domain $D\subseteq\CC$ the algebra
$\cO(\CC)$ may be viewed (via the restriction map)
as a holomorphically closed subalgebra of $\cO(D)$,
but it is not spectrally invariant in $\cO(D)$ unless $D=\CC$.
\end{remark}

It is clear from the definition that the intersection of any family of holomorphically
closed subalgebras of $A$ is holomorphically closed. This leads naturally to the
following definition.

\begin{definition}
The {\em holomorphic closure}, $\Hol(S)$, of a subset $S\subseteq A$ is
the intersection of all holomorphically closed subalgebras of $A$ containing $S$.
\end{definition}

Clearly, $\Hol(S)$ is the smallest holomorphically closed subalgebra of $A$ containing $S$.
It can also be described more explicitly as follows.

\begin{prop}
\label{prop:hol_closure}
For each subset $S\subseteq A$ we have
\[
\Hol(S)=\{ f(a) : f\in\cF _n,\; a\in S^n,\; n\in\Z_+\}.
\]
If, in addition, $A$ is commutative, then
\[
\Hol(S)=\{ f(a) : f\in\cO(\CC^n),\; a\in S^n,\; n\in\Z_+\}.
\]
\end{prop}

Now we are ready to introduce our main objects of study.

\begin{definition}
We say that $A$ is {\em holomorphically generated} by a subset $S\subseteq A$ if
$\Hol(S)=A$. We say that $A$ is {\em holomorphically finitely generated} (HFG for short)
if $A$ is holomorphically generated by a finite subset.
\end{definition}

Proposition~\ref{prop:hol_closure}, together with the Open Mapping Theorem,
implies the following characterization of HFG algebras.

\begin{prop}
\label{prop:HFG_quot_F}
Let $A$ be a Fr\'echet-Arens-Michael algebra.
\begin{mycompactenum}
\item $A$ is holomorphically finitely generated if and only if
$A$ is topologically isomorphic to $\cF_n/I$ for some $n\in\Z_+$ and for some closed
two-sided ideal $I\subseteq \cF_n$.
\item If $A$ is commutative, then $A$ is holomorphically finitely generated if and only if
$A$ is topologically isomorphic to $\cO(\CC^n)/I$ for some $n\in\Z_+$ and for some closed
two-sided ideal $I\subseteq \cO(\CC^n)$.
\end{mycompactenum}
\end{prop}

\begin{corollary}
Each HFG algebra is nuclear.
\end{corollary}

Combining Proposition~\ref{prop:HFG_quot_F} with the Remmert-Bishop-Narasimhan-Wiegmann
Embedding Theorem for Stein spaces \cite{Wiegmann}
and with Forster's theorem (see Example~\ref{ex:For_thm}), we obtain the following.

\begin{theorem}
\label{thm:comm_HFG}
A commutative Fr\'echet-Arens-Michael algebra is holomorphically finitely generated if and only if
it is topologically isomorphic to $\cO(X)$ for some Stein space $(X,\cO_X)$ of finite
embedding dimension. Moreover, the functor $(X,\cO_X)\mapsto\cO(X)$ is an anti-equivalence
between the category of Stein spaces of finite embedding dimension and the category
of commutative HFG algebras.
\end{theorem}

Theorem~\ref{thm:comm_HFG} looks similar to the Gelfand-Naimark Theorem
(Example~\ref{ex:GN}) and to the categorical consequence
of the Nullstellensatz (Example~\ref{ex:Null}). Thus HFG algebras may be considered
as candidates for ``noncommutative Stein spaces of finite embedding dimension''.
Of course, this naive point of view needs a solid justification, and it is perhaps too early to say whether
it will lead to an interesting theory. As a first step towards this goal, we will give some concrete
examples of HFG algebras below, showing that the class of HFG algebras is rather large.

\section{Arens-Michael envelopes}
\label{sect:AM_env}

Many natural examples of HFG algebras can be obtained from the following general construction.

\begin{definition}
Let $A$ be an algebra. A pair
$(\wh{A},i_A)$ consisting of an Arens-Michael algebra $\wh{A}$ and a homomorphism
$i_A\colon A\to\wh{A}$ is called the \emph{Arens-Michael envelope} of $A$ if for each
Arens-Michael algebra $B$ and for each homomorphism $\varphi\colon A\to B$
there exists a unique continuous homomorphism $\wh{\varphi}\colon\wh{A}\to B$
making the following diagram commute:
\[
\xymatrix{
\wh{A} \ar@{-->}[r]^{\wh{\varphi}} & B \\
A \ar[u]^{i_A} \ar[ur]_\varphi
}
\]
\end{definition}

Arens-Michael envelopes were introduced by J.~L.~Taylor \cite{T1} under the name of
``completed locally $m$-convex envelopes''. Now it is customary to call them
``Arens-Michael envelopes'', following the terminology suggested by A.~Ya.~Helemskii \cite{X2}.
For a more detailed study of Arens-Michael envelopes we refer to
\cite{Pir_qfree,Pir_stbflat,Pir_locsolv,Dos_holfun,Dos_Slodk,Dos_hdim,Dos_locleft,Dos_absbas}.

It is clear from the definition that if $\wh{A}$ exists, then it is unique up to a unique
topological isomorphism over $A$. In fact, $\wh{A}$ always exists, and it can be obtained
as the completion of $A$ with respect to the family of all submultiplicative seminorms on $A$.
Obviously, the correspondence $A\mapsto\wh{A}$ is a functor from the category of algebras
to the category of Arens-Michael algebras, and this functor is the left adjoint to the forgetful
functor acting in the opposite direction. In what follows, the functor $A\mapsto\wh{A}$
will be called the {\em Arens-Michael functor}.

Here are some basic examples of Arens-Michael envelopes.

\begin{example}[\cite{T2}]
\label{ex:AM_poly}
If $A=\CC[x_1,\ldots ,x_n]$ is the polynomial algebra, then $\wh{A}=\cO(\CC^n)$, the Fr\'echet
algebra of entire functions on $\CC^n$.
\end{example}

\begin{example}[\cite{Pir_qfree}]
\label{ex:AM_O}
If $X$ is an affine algebraic variety and $A$ is the algebra $\cO^\reg(X)$ of regular functions on $X$,
then $\wh{A}$ is the algebra $\cO^\hol(X)$ of holomorphic functions on $X$.
More generally, if $(X,\cO^\reg_X)$ is an affine scheme of finite type over $\CC$, and if
$A=\cO^\reg(X)$, then $\wh{A}=\cO^\hol(X)$, where
$(X,\cO^\hol_X)$ is the Stein space associated to $(X,\cO^\reg_X)$
(cf. \cite[Appendix~B]{Hart_AG}).
\end{example}

\begin{example}[\cite{T2}]
\label{ex:AM_F}
The Arens-Michael envelope of the free algebra $F_n$ is the algebra  $\cF_n$
of free entire functions \eqref{Fn}.
\end{example}

By combining Example~\ref{ex:AM_F} with Proposition~\ref{prop:HFG_quot_F}, and by using the fact that
the Arens-Michael functor commutes with quotients \cite{Pir_qfree}, we obtain the following.

\begin{corollary}
If $A$ is a finitely generated algebra, then $\wh{A}$ is holomorphically finitely generated.
\end{corollary}

Thus we may interpret the Arens-Michael functor as a ``noncommutative analytization functor''
acting from the category of noncommutative affine schemes of finite type to the category
of noncommutative Stein spaces of finite embedding dimension.
Example~\ref{ex:AM_O} shows that, in the commutative case, the Arens-Michael functor
reduces to the classical analytization functor $(X,\cO^\reg_X)\mapsto (X,\cO^\hol_X)$.

Here are some more examples of Arens-Michael envelopes.

\subsection{Quantum affine space}
Let $\CC^\times=\CC\setminus\{ 0\}$, and let
$q\in\CC^\times$. Recall that the algebra $\cO_q^\reg(\CC^n)$ {\em of regular functions on the
quantum affine $n$-space} is generated by
$n$ elements $x_1,\ldots ,x_n$ subject to the relations $x_i x_j=qx_j x_i$ for all $i<j$
(see, e.g., \cite{Br_Good}). If $q=1$, then $\cO_q^\reg(\CC^n)$ is nothing but the
polynomial algebra $\CC[x_1,\ldots ,x_n]=\cO^\reg(\CC^n)$.
Of course, $\cO_q^\reg(\CC^n)$ is noncommutative unless $q=1$, but the monomials
$x^\alpha=x_1^{\alpha_1}\cdots x_n^{\alpha_n}\; (\alpha\in\Z_+^n)$ still
form a basis of $\cO_q^\reg(\CC^n)$. Thus $\cO_q^\reg(\CC^n)$ may be viewed as
a ``deformed'' polynomial algebra.

Taking into account Example~\ref{ex:AM_poly}, it is natural to give the following definition.

\begin{definition}
The Arens-Michael envelope of $\cO_q^\reg(\CC^n)$
is denoted by $\cO_q^\hol(\CC^n)$ and is called the
{\em algebra of holomorphic functions on the quantum affine $n$-space}.
\end{definition}

The algebra $\cO_q^\hol(\CC^n)$ has the following explicit description.
Define a weight function $w_q\colon\Z_+^n\to\R_+$ by
\begin{equation}
\label{w_q}
w_q(\alpha)=
\begin{cases}
1 & \text{if } |q|\ge 1,\\
|q|^{\sum_{i<j}\alpha_i\alpha_j} & \text{if } |q|<1.
\end{cases}
\end{equation}

\begin{theorem}[\cite{Pir_qfree}]
\label{thm:q_aff_hol}
For each $q\in\CC^\times$ we have
\[
\cO^\hol_q(\CC^n)=
\Bigl\{
a=\sum_{\alpha\in\Z_+^n} c_\alpha x^\alpha :
\| a\|_\rho=\sum_{\alpha\in\Z_+^n} |c_\alpha| w_q(\alpha) \rho^{|\alpha|}<\infty
\;\forall\rho>0
\Bigr\}
\]
(where $|\alpha|=\alpha_1+\cdots+\alpha_n$ for $\alpha=(\alpha_1,\ldots ,\alpha_n)\in \Z_+^n$).
The topology on $\cO^\hol_q(\CC^n)$ is given by the seminorms $\|\cdot\|_\rho \; (\rho>0)$.
\end{theorem}

\subsection{Quantum torus}
Let $q\in\CC^\times$. Recall that the algebra $\cO_q^\reg((\CC^\times)^n)$
{\em of regular functions on the quantum $n$-torus} is generated by
$2n$ elements $x_1,\ldots ,x_n,x_1^{-1},\ldots ,x_n^{-1}$
subject to the relations $x_i^{-1} x_i=x_i x_i^{-1}=1$ for all $i$,
and $x_i x_j=qx_j x_i$ for all $i<j$
(see, e.g., \cite{Br_Good}). If $q=1$, then $\cO_q^\reg((\CC^\times)^n)$ is nothing but the
Laurent polynomial algebra $\CC[x_1^{\pm 1},\ldots ,x_n^{\pm 1}]=\cO^\reg((\CC^\times)^n)$.
Of course, $\cO_q^\reg((\CC^\times)^n)$ is noncommutative unless $q=1$, but the monomials
$x^\alpha=x_1^{\alpha_1}\cdots x_n^{\alpha_n}\; (\alpha\in\Z^n)$ still
form a basis of $\cO_q^\reg(\CC^n)$. Thus $\cO_q^\reg((\CC^\times)^n)$ may be viewed as
a ``deformed'' Laurent polynomial algebra.

\begin{definition}
If $|q|=1$, then the Arens-Michael envelope of $\cO_q^\reg((\CC^\times)^n)$
is denoted by $\cO_q^\hol((\CC^\times)^n)$ and is called the
{\em algebra of holomorphic functions on the quantum affine $n$-torus}.
\end{definition}

\begin{theorem}[\cite{Pir_qfree}]
\label{thm:q_tor_hol}
If $|q|=1$, then
\[
\cO^\hol_q((\CC^\times)^n)=
\Bigl\{
a=\sum_{\alpha\in\Z_+^n} c_\alpha x^\alpha :
\| a\|_\rho=\sum_{\alpha\in\Z^n} |c_\alpha| \rho^{|\alpha|}<\infty
\;\forall\rho>0
\Bigr\}.
\]
The topology on $\cO^\hol_q((\CC^\times)^n)$ is given by the seminorms $\|\cdot\|_\rho \; (\rho>0)$.
\end{theorem}

\begin{remark}
If $|q|\ne 1$, then the only submultiplicative seminorm on $\cO^\reg_q((\CC^\times)^n)$
is identically zero \cite{Pir_qfree}, and so the Arens-Michael envelope of
$\cO^\reg_q((\CC^\times)^n)$ is trivial. Another example of an algebra with trivial
Arens-Michael envelope is the Weyl algebra generated by two elements $x,y$
subject to the relation $[x,y]=1$.
\end{remark}

We refer to \cite{Pir_qfree,Dos_hdim,Dos_absbas} for explicit descriptions of Arens-Michael envelopes
of some other algebras, including quantum Weyl algebras,
the algebra of quantum $2\times 2$-matrices, and universal enveloping algebras.

\section{Free products and free polydisk}

To construct more examples of HFG algebras (not necessarily related
to Arens-Michael envelopes), we need the simple fact that the category of HFG algebras
is stable under the analytic free product.

\begin{definition}[cf. \cite{Cuntz_doc}]
Let $A_1$ and $A_2$ be Arens-Michael algebras. The {\em analytic free product} of $A_1$ and $A_2$
is the coproduct of $A_1$ and $A_2$ in the category of Arens-Michael algebras.
More explicitly, the analytic free product of $A_1$ and $A_2$ is a triple
$(A_1\Pfree A_2,j_1,j_2)$
consisting of an Arens-Michael algebra $A_1\Pfree A_2$ and continuous homomorphisms
$j_k\colon A_k\to A_1\Pfree A_2\; (k=1,2)$ such that for each Arens-Michael algebra
$B$ and each pair $\varphi_k\colon A_k\to B\; (k=1,2)$ of continuous homomorphisms
there exists a unique continuous homomorphism
$A_1\Pfree A_2\to B$ making the following diagram commute:
\[
\xymatrix@-10pt{
& A_1\Pfree A_2 \ar@{-->}[dd] \\
A_1 \ar[ur]^{j_1} \ar[dr]_{\varphi_1} && A_2 \ar[ul]_{j_2} \ar[dl]^{\varphi_2} \\
& B
}
\]
\end{definition}

Clearly, if $A_1\Pfree A_2$ exists, then it is unique up to a unique
topological algebra isomorphism over $A_1$ and $A_2$.
To show that $A_1\Pfree A_2$ exists, recall the definition
(due to J.~Cuntz~\cite{Cuntz_doc}; cf. also \cite{Pir_qfree})
of an analytic tensor algebra.
Let $E$ be a complete locally convex space, and let $\{\|\cdot\|_\lambda : \lambda\in\Lambda\}$
be a directed family of seminorms generating the topology on $E$.
For each $\lambda\in\Lambda$ and each $n\in\Z_+$, let $\|\cdot\|_\lambda^{(n)}$
denote the $n$th projective tensor power of $\|\cdot\|_\lambda$
(we let $\|\cdot\|_\lambda^{(0)}=|\cdot|$ by definition).
The {\em analytic tensor algebra} $\wh{T}(E)$ is defined by
\[
\wh{T}(E)=\Bigl\{ a=\sum_{n=0}^\infty a_n : a_n\in E^{\wh{\otimes}n},\;
\| a\|_{\lambda,\rho}=\sum_n \| a_n\|_\lambda^{(n)} \rho^n<\infty\; \forall\lambda\in\Lambda,
\;\forall\rho>0\Bigr\}.
\]
The topology on $\wh{T}(E)$ is defined by the seminorms
$\|\cdot\|_{\lambda,\rho}\; (\lambda\in\Lambda,\;\rho>0)$, and the product
on $\wh{T}(E)$ is given by concatenation, like on the usual tensor algebra $T(E)$.
Each seminorm $\|\cdot\|_{\lambda,\rho}$ is easily seen to be
submultiplicative, and so $\wh{T}(E)$ is an Arens-Michael algebra containing
$T(E)$ as a dense subalgebra. As was observed by Cuntz \cite{Cuntz_doc},
$\wh{T}(E)$ has the universal property that, for every Arens-Michael algebra $A$,
each continuous linear map $E\to A$ uniquely extends to a continuous
homomorphism $\wh{T}(E)\to A$. Note that $\wh{T}(\CC^n)\cong\cF_n$,
and that the free entire functional calculus \eqref{free_entire}
is a special case of the universal property of $\wh{T}(E)$.

Now let $A_1\Pfree A_2$ be the completion of the quotient of $\wh{T}=\wh{T}(A_1\oplus A_2)$
by the two-sided closed ideal generated by all elements of the form
\[
a_i\otimes b_i-a_i b_i,\quad 1_{\wh{T}}-1_{A_i}\quad
(a_i,b_i\in A_i,\; i=1,2).
\]
The canonical homomorphisms $j_k\colon A_k\to A_1\Pfree A_2$ are defined in
the obvious way. Using the universal property of $\wh{T}$, it is easy to show that
$A_1\Pfree A_2$ is indeed the analytic free product of $A_1$ and $A_2$.

A more explicit construction
of $A_1\Pfree A_2$ (in the nonunital category) was given by Cuntz~\cite{Cuntz_doc}.
We adapt the construction to the unital case below, assuming that for each $i=1,2$
there exists a closed two-sided ideal $\ol{A}_i\subset A_i$ such that
$A_i=\ol{A}_i\oplus\CC 1_{A_i}$ as locally convex spaces.

For each $k\ge 2$, let
\[
W'_{2,k}=\Bigl\{ \alpha=(\alpha_1,\ldots ,\alpha_k)\in W_{2,k} : \alpha_i\ne\alpha_{i+1}
\;\forall i=1,\ldots ,k-1\Bigr\}.
\]
Let also $W'_{2,1}=W_{2,1}$, $W'_{2,0}=W_{2,0}$, and $W'_2=\bigsqcup_{k\in\Z_+} W'_{2,k}$.
Given $\alpha=(\alpha_1,\ldots ,\alpha_k)\in W'_2$, let
$\ol{A}_\alpha=\ol{A}_{\alpha_1}\Ptens\cdots\Ptens \ol{A}_{\alpha_k}$.
For $k=0$, we let $\ol{A}_\emptyset=\CC$.
Choose directed defining families $\{\|\cdot\|_{\lambda_i} : \lambda_i\in\Lambda_i\}\; (i=1,2)$
of submultiplicative seminorms on $\ol{A}_1$ and $\ol{A}_2$, respectively.
Let $\Lambda=\Lambda_1\times\Lambda_2$, and, for each
$\lambda=(\lambda_1,\lambda_2)\in\Lambda$ and each $\alpha=(\alpha_1,\ldots ,\alpha_k)\in W'_2$,
let $\|\cdot\|_\lambda^{(\alpha)}$ denote the projective tensor product of the
seminorms $\|\cdot\|_{\lambda_{\alpha_1}},\ldots ,\|\cdot\|_{\lambda_{\alpha_k}}$.
For $k=0$, it is convenient to set $\|\cdot\|_\lambda^{(\emptyset)}=|\cdot|$.
We have
\begin{equation}
\label{Pfree_expl}
A_1\Pfree A_2=
\Bigl\{ a=\sum_{\alpha\in W'_2} a_\alpha : a_\alpha\in \ol{A}_\alpha,\;
\| a\|_{\lambda,\rho}=\sum_\alpha \| a_\alpha\|_\lambda^{(\alpha)} \rho^{|\alpha|}<\infty
\;\forall \lambda\in\Lambda,\; \forall\rho>0\Bigr\}.
\end{equation}
The topology on $A_1\Pfree A_2$ is defined by the seminorms
$\|\cdot\|_{\lambda,\rho}\; (\lambda\in\Lambda,\rho>0)$.
The product on $A_1\Pfree A_2$ is given by concatenation composed (if necessary) with the
product maps $\ol{A}_i\Ptens \ol{A}_i\to \ol{A}_i\; (i=1,2)$.

\begin{example}
\label{ex:Fn-Fm}
If $E$ and $F$ are complete locally convex spaces, then we have a topological
algebra isomorphism $\wh{T}(E)\Pfree\wh{T}(F)\cong\wh{T}(E\oplus F)$.
In particular, $\cF_m\Pfree\cF_n\cong\cF_{m+n}$, and
\begin{equation}
\label{F_O_free}
\cF_n\cong\cO(\CC)\Pfree\cdots\Pfree\cO(\CC).
\end{equation}
\end{example}

The next result easily follows from Example~\ref{ex:Fn-Fm} and
Proposition~\ref{prop:HFG_quot_F}.

\begin{prop}
\label{prop:HFG_free}
If $A_1$ and $A_2$ are HFG algebras, then so is $A_1\Pfree A_2$.
\end{prop}

Using Proposition~\ref{prop:HFG_free}, we can construct more examples of HFG
algebras. Let $r>0$, and let $\DD_r=\{ z\in\CC : |z|<r\}$ denote the open disk of radius $r$.
The following definition is motivated by~\eqref{F_O_free}.

\begin{definition}
We define the {\em algebra of holomorphic functions on the free $n$-dimensional polydisk of
radius $r$} to be
\begin{equation}
\label{F_poly}
\cF(\DD^n_r)=\cO(\DD_r)\Pfree\cdots\Pfree\cO(\DD_r).
\end{equation}
\end{definition}

By Proposition~\ref{prop:HFG_free}, $\cF(\DD^n_r)$ is an HFG algebra.
Note that replacing in \eqref{F_O_free} and \eqref{F_poly} the analytic free product $\Pfree$
by the projective tensor product $\Ptens$ yields the algebras of holomorphic functions
on $\CC^n$ and $\DD_r^n$, respectively.

We have a canonical ``restriction'' map $\cF_n\to\cF(\DD^n_r)$ defined to be the
$n$th free power of the restriction map $\cO(\CC)\to\cO(\DD_r)$.
For each $i=1,\ldots ,n$, the canonical image of the
free generator $\zeta_i\in\cF_n$ in $\cF(\DD^n_r)$ will also be denoted by $\zeta_i$.
It is easy to see that $\cF(\DD^n_r)$ has the following universal property.

\begin{prop}
\label{prop:univ_F_poly}
Let $A$ be an Arens-Michael algebra, and let $a=(a_1,\ldots ,a_n)$ be an $n$-tuple in $A^n$
such that $\sigma_A(a_i)\subseteq\DD_r$ for all $i=1,\ldots ,n$. Then there exists a unique
continuous homomorphism $\gamma_a^\free\colon\cF(\DD^n_r)\to A$
such that $\gamma_a^\free(\zeta_i)=a_i$ for all $i=1,\ldots ,n$.
\end{prop}

The algebra $\cF(\DD_r^n)$ can also be described more explicitly as follows.
Given $k\ge 2$ and $\alpha=(\alpha_1,\ldots ,\alpha_k)\in W_n$, let $d(\alpha)$ denote the
cardinality of the set
\[
\bigl\{ i \in \{ 1,\ldots ,k-1\} : \alpha_i\ne\alpha_{i+1} \bigr\}.
\]
If $|\alpha|\in\{ 0,1\}$, we set $d(\alpha)=|\alpha|-1$.
The next result follows from \eqref{Pfree_expl}.

\begin{prop}
We have
\[
\cF(\DD_r^n)=\Bigl\{ a=\sum_{\alpha\in W_n} c_\alpha\zeta_\alpha :
\| a\|_{\rho_1,\rho_2}=\sum_{\alpha\in W_n} |c_\alpha|\rho_1^{|\alpha|} \rho_2^{d(\alpha)+1}<\infty
\;\forall 0<\rho_1<r,\; \forall \rho_2>0\Bigr\}.
\]
The topology on $\cF(\DD_r^n)$ is given by the seminorms
$\|\cdot\|_{\rho_1,\rho_2}\; (0<\rho_1<r,\; \rho_2>0)$,
and the product is given by concatenation.
\end{prop}

\begin{remark}
\label{rem:Tay_alg}
It is instructive to compare $\cF(\DD_r^n)$ with Taylor's
``free power series algebras'' $\cF_n(r)$ \cite{T2,T3}.
By definition,
\[
\cF_n(r)=\Bigl\{ a=\sum_{\alpha\in W_n} c_\alpha\zeta_\alpha :
\| a\|_\rho=\sum_{\alpha\in W_n} |c_\alpha|\rho^{|\alpha|}<\infty
\;\forall 0<\rho<r\Bigr\}.
\]
Comparing \eqref{Fn}, \eqref{F_O_free}, and \eqref{F_poly}, we see that
$\cF_n(r)=\cF_n=\cF(\DD_r^n)$ if $r=\infty$. However, if $r<\infty$, then
$\cF(\DD_r^n)$ differs from $\cF_n(r)$. Indeed, $\cF_n(r)$ is not nuclear \cite{Lum},
and so it is not even an HFG algebra. Although such algebras fall outside the
scope of this paper, we would like to make a short digression and to show
that $\cF_n(r)$ has a remarkable
universal property, similar in spirit to Proposition~\ref{prop:univ_F_poly}.

Recall that, for a Banach algebra $A$ and an $n$-tuple $a=(a_1,\ldots ,a_n)\in A^n$,
the {\em joint spectral radius} $r_\infty(a)$ is given by
\[
r_\infty(a)=\limsup_{k\to\infty}\Bigl(\sup_{\alpha\in W_{n,k}} \| a_\alpha\|\Bigr)^{1/k}.
\]
If now $A$ is an Arens-Michael algebra, we say that an $n$-tuple $a\in A^n$ is
{\em strictly $r$-contractive} if, for each Banach algebra $B$ and each continuous
homomorphism $\varphi\colon A\to B$, we have $r_\infty(\varphi^{\times n}(a))<r$.
An equivalent, but more handy definition is as follows. Let $\{\|\cdot\|_\lambda : \lambda\in\Lambda\}$
be a directed defining family of submultiplicative seminorms on $A$. For each $\lambda\in\Lambda$,
let $A_\lambda$ denote the completion of $A$ with respect to $\|\cdot\|_\lambda$, and let
$a_\lambda$ denote the canonical image of $a$ in $A_\lambda^n$.
Then it is easy to show that $a$ is strictly $r$-contractive if and only if $r_\infty(a_\lambda)<r$
for all $\lambda\in\Lambda$. For example, if $\zeta=(\zeta_1,\ldots ,\zeta_n)\in\cF_n(r)^n$
is the canonical system of free generators, then $\zeta$ is strictly $r$-contractive, but is
not strictly $\rho$-contractive for $\rho<r$.

Now it is easy to show that $\cF_n(r)$ has the following universal property.
{\em Given an Arens-Michael algebra $A$ and a strictly $r$-contractive
$n$-tuple $a=(a_1,\ldots ,a_n)\in A^n$, there exists a unique
continuous homomorphism $\gamma_a^\free\colon\cF_n(r)\to A$
such that $\gamma_a^\free(\zeta_i)=a_i$ for all $i=1,\ldots ,n$.}

To complete our discussion of $\cF_n(r)$ and $\cF(\DD_r^n)$, let us note that
the $n$-tuple $(\zeta_1,\ldots ,\zeta_n)$ is not strictly $R$-contractive
in $\cF(\DD_r^n)^n$ for any $R>0$.
\end{remark}

\section{Quantum polydisk and quantum ball}

\begin{definition}
Let $q\in\CC^\times$, and let $r>0$.
We define the {\em algebra of holomorphic functions on the
quantum $n$-polydisk of radius $r$} by
\[
\cO^\hol_q(\DD^n_r)=
\Bigl\{
a=\sum_{\alpha\in\Z_+^n} c_\alpha x^\alpha :
\| a\|_\rho=\sum_{\alpha\in\Z_+^n} |c_\alpha| w_q(\alpha) \rho^{|\alpha|}<\infty
\;\forall 0<\rho<r
\Bigr\},
\]
where the function $w_q\colon \Z_+^n\to\R_+$ is given by \eqref{w_q}.
The topology on $\cO^\hol_q(\DD^n_r)$ is defined by the seminorms $\|\cdot\|_\rho\; (0<\rho<r)$,
and the product on $\cO^\hol_q(\DD^n_r)$ is uniquely determined by $x_i x_j=qx_j x_i$ for all $i<j$.
\end{definition}

In other words, $\cO^\hol_q(\DD^n_r)$ is the completion of $\cO_q^\reg(\CC^n)$
(see Section~\ref{sect:AM_env}) with respect
to the family $\{\|\cdot\|_\rho : 0<\rho<r\}$ of submultiplicative seminorms.
Clearly, if $q=1$, then $\cO^\hol_q(\DD^n_r)$ is topologically isomorphic to the algebra
$\cO^\hol(\DD^n_r)$ of holomorphic functions on the polydisk $\DD_r^n$.

\begin{prop}
\label{prop:q_poly_quot_free_poly}
For each $q\in\CC^\times$, $\cO^\hol_q(\DD^n_r)$ is topologically isomorphic to the
quotient of $\cF(\DD_r^n)$ by the two-sided closed ideal generated by the elements
$\zeta_i\zeta_j-q\zeta_j\zeta_i$ for all $i<j$. As a corollary, $\cO^\hol_q(\DD^n_r)$
is an HFG algebra.
\end{prop}

Note that $\cF(\DD_r^n)$ can be replaced by $\cF_n(r)$ in the above proposition.

Our last example is the algebra of holomorphic functions on the quantum ball.
The construction below is a slight modification of L.~L.~Vaksman's
$q$-analogue of $A(\bar\BB^n)$, the algebra of functions continuous on
the closed ball and holomorphic on the open ball~\cite{Vaks_max}.
In order to motivate the construction, let us start from the classical situation.
Let $z_1,\ldots ,z_n$ be the complex coordinates on $\CC^n$, and let
$\BB^n_r=\{ z\in\CC^n : \sum_i |z_i|^2<r^2\}$ denote the open ball
of radius $r$. For brevity, we write $\BB^n$ for $\BB^n_1$.
Denote by $\Pol(\CC^n)$ the algebra of polynomials in $z_1,\ldots ,z_n$ and
their complex conjugates $\bar z_1,\ldots ,\bar z_n$. There is a natural involution
on $\Pol(\CC^n)$ uniquely determined by $z_i^*=\bar z_i\; (i=1,\ldots ,n)$.
By the Stone-Weierstrass Theorem, the completion of $\Pol(\CC^n)$ with respect
to the uniform norm $\| f\|=\sup_{z\in\bar\BB^n} |f(z)|$ is the algebra
$C(\bar\BB^n)$ of continuous functions on the closed ball $\bar\BB^n$.
For each $r\in (0,1)$ and each $f\in \CC[z_1,\ldots ,z_n]$, let
\[
\| f\|_r=\sup_{z\in\bar\BB^n_r} |f(z)| = \|\gamma_r(f)\|,
\]
where $\gamma_r$ is the automorphism of $\CC[z_1,\ldots ,z_n]$ uniquely
determined by $\gamma_r(z_i)=rz_i\; (i=1,\ldots ,n)$.
It is easy to see that the completion of $\CC[z_1,\ldots ,z_n]$ with respect to
the family $\{\|\cdot\|_r : 0<r<1\}$ of seminorms is topologically isomorphic
to the algebra $\cO(\BB^n)$ of holomorphic functions on $\BB^n$.

Now let us ``quantize'' the above data.
Fix $q\in (0,1)$, and denote by $\Pol_q(\CC^n)$
the $*$-algebra generated (as a $*$-algebra) by $n$ elements $z_1,\ldots ,z_n$
subject to the relations
\begin{equation}
\label{tw_CCR}
\begin{aligned}
z_i z_j&=q z_j z_i \quad (i<j);\\
z_i^* z_j&= q z_j z_i^* \quad (i\ne j);\\
z_i^* z_i&=q^2 z_i z_i^*+(1-q^2)\Bigl(1-\sum_{k>i} z_k z_k^*\Bigr).
\end{aligned}
\end{equation}
Clearly, for $q=1$ we have $\Pol_q(\CC^n)=\Pol(\CC^n)$.
The algebra $\Pol_q(\CC^n)$ was introduced by W.~Pusz and S.~L.~Woronowicz \cite{PW},
although they used different $*$-generators $a_1,\ldots ,a_n$ given by
$a_i=(1-q^2)^{-1/2} z_i^*$. Relations \eqref{tw_CCR} divided by $1-q^2$ and
written in terms of the $a_i$'s
are called the ``twisted canonical commutation relations'', and the algebra $A_q=\Pol_q(\CC^n)$
defined in terms of the $a_i$'s is sometimes called the ``quantum Weyl algebra''
(see, e.g., \cite{WZ,Alev,Jordan,Klim_Schm}).
Note that, while $\Pol_q(\CC^n)$ becomes $\Pol(\CC^n)$ for $q=1$, $A_q$ becomes
the Weyl algebra. The idea to use the generators $z_i$
instead of the $a_i$'s and to consider $\Pol_q(\CC^n)$ as a $q$-analogue of $\Pol(\CC^n)$
is probably due to Vaksman \cite{Vaks_splet}; the one-dimensional case
was considered in \cite{Klim_Lesn}.
The algebra $\Pol_q(\CC^n)$ serves as a basic example in the general theory
of ``quantum bounded symmetric domains'' developed by Vaksman and his
collaborators (see \cite{Vaks_sborn,Vaks_book} and references therein).

The $q$-analogue of the algebra $C(\bar\BB^n)$ is defined as follows.
For each $k\in\N$, let $[k]_q=\sum_{i=0}^{k-1} q^{2i}$.
Fix a separable Hilbert space $H$ with an orthonormal basis
$\{ e_\alpha : \alpha\in \Z_+^n\}$. Following \cite{PW}, for each
$\alpha=(\alpha_1,\ldots ,\alpha_n)\in\Z_+^n$ we will write $|\alpha_1,\ldots ,\alpha_n\ra$
for $e_\alpha$. As was proved by Pusz and Woronowicz \cite{PW}, there exists
a faithful irreducible $*$-representation $\pi\colon\Pol_q(\CC^n)\to\cB(H)$ uniquely
determined by
\begin{gather*}
\pi(z_j)e_\alpha=\sqrt{1-q^2} \sqrt{[\alpha_j+1]} q^{\sum_{k>j}\alpha_k}
|\alpha_1,\ldots ,\alpha_j+1,\ldots ,\alpha_n\ra\\
(j=1,\ldots ,n,\; \alpha=(\alpha_1,\ldots ,\alpha_n)\in\Z_+^n).
\end{gather*}
The completion of $\Pol_q(\CC^n)$ with respect to the operator norm $\| a\|_\op=\| \pi(a)\|$
is denoted by $C_q(\bar\BB^n)$ and is called the {\em algebra of continuous functions
on the closed quantum ball} \cite{Vaks_max}; see also \cite{PW,Prosk_Sam}.

Now we are ready to define the algebra of holomorphic functions on the quantum ball.
Observe that the subalgebra of $\Pol_q(\CC^n)$ generated by $z_1,\ldots ,z_n$
is exactly $\cO_q^\reg(\CC^n)$. For each $r\in (0,1)$, let $\gamma_r$ be the automorphism
of $\cO_q^\reg(\CC^n)$ uniquely
determined by $\gamma_r(z_i)=rz_i\; (i=1,\ldots ,n)$. Define a norm $\|\cdot\|_r$ on
$\cO_q^\reg(\CC^n)$ by $\| a\|_r=\|\gamma_r(a)\|_\op$.

\begin{definition}
The completion of $\cO_q^\reg(\CC^n)$ with respect to the family $\{\|\cdot\|_r : 0<r<1\}$
of seminorms is denoted by $\cO^\hol_q(\BB^n)$ and is called the {\em algebra of holomorphic
functions on the quantum ball}.
\end{definition}

It follows from the above discussion that, if we replace
$C_q(\bar\BB^n)$ by $C(\bar\BB^n)$ in the above construction, then
the result will be the algebra $\cO^\hol(\BB^n)$ of holomorphic functions on $\BB^n$.
In Section~\ref{sect:def}, we will study relations between $\cO_q^\hol(\BB^n)$ and $\cO^\hol(\BB^n)$
in more detail.

At the moment, it is not obvious whether $\cO_q^\hol(\BB^n)$ is an HFG algebra.
The positive answer follows from the next result. Let $\DD^n=\DD^n_1$ denote the open
unit polydisk in $\CC^n$.

\begin{theorem}
\label{thm:poly_vs_ball}
For each $0<q<1$
the algebras $\cO_q^\hol(\BB^n)$ and $\cO_q^\hol(\DD^n)$ are topologically isomorphic.
\end{theorem}

Of course, Theorem~\ref{thm:poly_vs_ball} is a ``purely quantum'' phenomenon.
Indeed, the classical function algebras $\cO^\hol(\BB^n)$ and $\cO^\hol(\DD^n)$ are not isomorphic
(unless $n=1$) because $\BB^n$ and $\DD^n$ are not biholomorphically equivalent.

Let us briefly discuss the idea of the proof.
Consider the following families of seminorms on $\cO_q^\reg(\CC^n)$:
\begin{align}
\label{semi_1}
\| a\|_r^{(1)}&=\sum_{\alpha\in\Z_+^n} |c_\alpha| w_q(\alpha) r^{|\alpha|}
\quad (0<r<1);\\
\label{semi_2}
\| a\|_r^{(2)}&=\Bigl(\sum_{\alpha\in\Z_+^n} |c_\alpha|^2 w_q^2(\alpha) r^{2|\alpha|}\Bigr)^{1/2}
\quad (0<r<1);\\
\label{semi_inf}
\| a\|_r^{(\infty)}&=\sup_{\alpha\in\Z_+^n} |c_\alpha| w_q(\alpha) r^{|\alpha|}
\quad (0<r<1);\\
\label{semi}
\| a\|_r&=\| \gamma_r(a)\|_\op \quad (0<r<1).
\end{align}
Recall that the completion of $\cO_q^\reg(\CC^n)$ w.r.t. the seminorms \eqref{semi_1}
is $\cO_q^\hol(\DD^n)$, while the completion of $\cO_q^\reg(\CC^n)$ w.r.t. the seminorms
\eqref{semi} is $\cO_q^\hol(\BB^n)$. A standard argument shows that the families
\eqref{semi_1}--\eqref{semi_inf} are equivalent. The following lemma gives a less
obvious estimate.

\begin{lemma}
\label{lemma:key}
We have
\begin{equation}
\label{key_est}
\Bigl(\prod_{j=1}^\infty (1-q^{2j})\Bigr)^{\frac{n}{2}} \|\cdot\|_1^{(2)}
\le \|\cdot\|_\op \le \|\cdot\|_1^{(1)}.
\end{equation}
\end{lemma}

The first inequality in \eqref{key_est} follows from the estimate $\| a\|_\op\ge \|\pi(a)e_0\|$,
where $e_0=|0,\ldots ,0\ra\in H$ is the ``vacuum vector''. The second inequality in \eqref{key_est}
follows from the fact that $\|\cdot\|_1^{(1)}$ is the largest submultiplicative seminorm on
$\cO_q^\reg(\CC^n)$ such that $\| z_i\|_1^{(1)}=1$ for all $i=1,\ldots ,n$
\cite[5.10]{Pir_qfree}.

Combining Lemma~\ref{lemma:key} with the equivalence of the families \eqref{semi_1} and
\eqref{semi_2}, we obtain the following.

\begin{corollary}
For each $0<\rho<r<1$ we have
\begin{equation}
\label{key2}
\Bigl(\frac{r^2-\rho^2}{r^2}\prod_{j=1}^\infty (1-q^{2j})\Bigr)^{\frac{n}{2}} \| \cdot\|_\rho^{(1)}
\le \| \cdot\|_r \le \| \cdot\|_r^{(1)}.
\end{equation}
\end{corollary}

Thus the families \eqref{semi_1} and \eqref{semi} of seminorms are equivalent, and so
$\cO_q^\hol(\BB^n)\cong\cO_q^\hol(\DD^n)$, as required.

Note that, while the second inequality in \eqref{key2} holds in the classical case $q=1$ as well,
the first inequality in \eqref{key2} has no classical counterpart.
Indeed, if we fix $r$ and take $\rho<r$ close enough to $r$, then the polydisk of
radius $\rho$ will not be contained in the ball of radius $r$, and so the supremum over the
polydisk will not be dominated by the supremum over the ball.

\section{Relation to Fr\'echet algebra bundles}
\label{sect:def}

In this section we discuss in which sense $\cO_q^\hol(\BB^n)$ and $\cO_q^\hol(\DD^n)$
are deformations of the function algebras $\cO^\hol(\BB^n)$ and $\cO^\hol(\DD^n)$,
respectively.

Recall some definitions from \cite{Gierz}. Let $X$ and $E$ be sets, and let
$p\colon E\to X$ be a surjective map. Suppose that for each $x\in X$ the
fiber $E_x=p^{-1}(x)$ is endowed with the structure of a vector space.
A function $\|\cdot\|\colon E\to [0,+\infty)$
is a {\em seminorm} if its restriction to each $E_x$ is a seminorm
in the usual sense. Let
\[
E\times_X E=\{ (x,y)\in E\times E : p(x)=p(y)\}.
\]
Suppose now that $X$ is a locally compact, Hausdorff topological space,
$A$ is a topological space, and $p\colon A\to X$ is a continuous, open surjection.
Suppose also that each fiber $A_x=p^{-1}(x)$ is endowed with the structure of an algebra
in such a way that the operations
\begin{align*}
A\times_X A &\to A,\quad (a,b)\mapsto a+b,\\
\CC\times A &\to A,\quad (\lambda,a)\mapsto\lambda a,\\
A\times_X A &\to A,\quad (a,b)\mapsto ab,
\end{align*}
are continuous. Let $\cN=\{ \|\cdot\|_i : i\in I\}$
be a directed family of seminorms on $A$ having an at most countable cofinal
subfamily. Assume that for each $x\in X$ the sets
\[
\{ a\in A : p(a)\in U,\; \| a\|_i<\eps\}\qquad
(i\in I,\; \eps>0,\; U\subseteq X\text{ is a neighbourhood of }x)
\]
form a neighbourhood base of $0\in A_x$
(this implies, in particular, that the topology on $A_x$ inherited from $A$
is determined by the seminorms $\|\cdot\|_i$, and that each seminorm $\|\cdot\|_i$
is an upper semicontinuous function on $A$).
Finally, assume that each fiber
$A_x$ is complete (so, according to the above remarks, $A_x$ is a Fr\'echet algebra).
Under the above assumptions, the pair $(A,p)$ together with the family $\cN$
is called a {\em Fr\'echet algebra bundle} over $X$.

\begin{theorem}
\label{thm:bnd}
\begin{mycompactenum}
\item For each $r\in (0,+\infty]$ there exists a Fr\'echet algebra bundle
$(D,p)$ over $\CC^\times$ such that for each $q\in\CC^\times$ we have
$D_q\cong\cO^\hol_q(\DD_r^n)$.
\item There exists a Fr\'echet algebra bundle
$(B,p)$ over $\CC^\times$ such that $B_1\cong \cO^\hol(\BB^n)$
and such that for each $q\in (0,1)$ we have
$B_q\cong\cO^\hol_q(\BB^n)$.
\end{mycompactenum}
\end{theorem}

We will not give a detailed proof here for lack of space.
Let us only explain how the bundles $(D,p)$ and $(B,p)$ are constructed.
By a {\em Fr\'echet $\cO(\CC^\times)$-algebra}
we mean a Fr\'echet algebra $R$ together with the structure of a left Fr\'echet
$\cO(\CC^\times)$-module such that the multiplication $R\times R\to R$ is $\cO(\CC^\times)$-bilinear.
For each $q\in\CC^\times$, let $M_q=\{ f\in\cO(\CC^\times) : f(q)=0\}$.
It is easy to see that, if $R$ is a Fr\'echet $\cO(\CC^\times)$-algebra, then $\ol{M_q\cdot R}$
is a closed two-sided ideal of $R$. The Fr\'echet algebra $R_q=R/\ol{M_q\cdot R}$
is called the {\em fiber} of $R$ over $q$. The following lemma explains the
terminology.

\begin{lemma}
Let $R$ be a Fr\'echet $\cO(\CC^\times)$-algebra, and let $A=\bigsqcup_{q\in\CC^\times} R_q$.
Define the projection $p\colon A\to\CC^\times$ by $p(a)=q$ if $a\in R_q$.
Let also $\{ \|\cdot\|_i^R : i\in I\}$ be a directed defining family of seminorms on $R$.
For each $q\in\CC^\times$ and each $i\in I$, let $\|\cdot\|_{i,q}$ denote the corresponding
quotient seminorm on $R_q$, and let $\|\cdot\|_i$ be the seminorm on $A$ whose restriction
to $R_q$ is $\|\cdot\|_{i,q}$. Then there is a unique topology on $A$ making
$(A,p)$ into a Fr\'echet algebra
bundle over $\CC^\times$ with respect to the family $\{\|\cdot\|_i : i\in I\}$ of
seminorms.
\end{lemma}

Now, in order to prove part (i) (respectively, (ii)) of Theorem~\ref{thm:bnd},
we need to construct a Fr\'echet $\cO(\CC^\times)$-algebra $R$ whose fiber over each $q\in\CC^\times$
(respectively, over each $q\in (0,1]$)
is isomorphic to $\cO_q^\hol(\DD^n_r)$ (respectively, to $\cO_q^\hol(\BB^n)$).
This can be done as follows. Let $z$ denote the complex coordinate on $\CC^\times$,
and let $I$ be the closed two-sided ideal of $\cO(\CC^\times,\cF(\DD^n_r))$
generated by the elements $x_i x_j-zx_j x_i\; (i<j)$. The quotient $R=\cO(\CC^\times,\cF(\DD^n_r))/I$
is a Fr\'echet $\cO(\CC^\times)$-algebra in a natural way, and
Proposition~\ref{prop:q_poly_quot_free_poly} easily implies that
$R_q\cong \cO_q^\hol(\DD^n_r)$ for each $q\in\CC^\times$. We can also use
Taylor's $\cF_n(r)$ instead of $\cF(\DD^n_r)$; the algebra $R$ will then be the same.

To prove part (ii), we have to replace $\cF(\DD^n_r)$ by
the algebra $\cF(\BB^n)$ of ``holomorphic functions on the free ball'' introduced by
G.~Popescu \cite{Pop_holball}\footnote[1]{Note that our $\cF(\BB^n)$ is Popescu's
$Hol(B(\mathcal H)_1^n)$; we have changed the notation in order to emphasize the
similarity between Popescu's algebras, our $\cF(\DD^n_r)$, and Taylor's $\cF_n(r)$.}. By definition,
\[
\cF(\BB^n)=
\Bigl\{ a=\sum_{\alpha\in W_n} c_\alpha\zeta_\alpha :
\limsup_{k\to\infty} \Bigl(\sum_{|\alpha|=k} |c_\alpha|^2\Bigr)^{1/2k}\le 1\Bigr\}.
\]
As was observed by Popescu \cite{Pop_holball},
$\cF(\BB^n)$ is indeed an algebra with respect to
the concatenation product.

Similarly to $\cF(\DD^n_r)$ and $\cF_n(r)$, the algebra
$\cF(\BB^n)$ can be used to construct a kind of ``free holomorphic functional
calculus'' (cf. Proposition~\ref{prop:univ_F_poly} and Remark~\ref{rem:Tay_alg}).
Let $H$ be a Hilbert space, and
let $T=(T_1,\ldots ,T_n)$ be an $n$-tuple in $\cB(H)^n$.
Following \cite{Pop_holball}, we identify $T$ with the ``row'' operator
acting from the Hilbert direct sum $H^n=H\oplus\cdots\oplus H$ to $H$.
Thus we have $\| T\|=\|\sum_{i=1}^n T_i T_i^*\|^{1/2}$.
If $\| T\|<1$, then for each $f=\sum_\alpha c_\alpha\zeta_\alpha\in\cF(\BB^n)$
the series $\sum_\alpha c_\alpha T_\alpha$ converges in $\cB(H)$
to an operator $f(T)$ \cite{Pop_holball}, and we have an algebra
homomorphism $\gamma_T^\free\colon\cF(\BB^n)\to\cB(H)$, $f\mapsto f(T)$.
In fact, $\gamma_T^\free$ exists under the weaker assumption
that Bunce's ``hilbertian'' joint spectral
radius \cite{Bunce_models} of $T$ is less
than $1$ (see \cite[4.1]{Pop_holball}), but we will not use this fact in the sequel.

There is a canonical topology on $\cF(\BB^n)$ defined as follows.
Fix an infinite-dimensional Hilbert space $H$, and, for each $r\in (0,1)$ and each
$f\in\cF(\BB^n)$, let
\[
\| f\|_r^{(P)}=\sup\bigl\{ \| f(T)\| : T\in\cB(H)^n,\; \| T\|\le r\bigr\}.
\]
By \cite[5.6]{Pop_holball}, $\cF(\BB^n)$ is a Fr\'echet space with respect
to the family $\{\|\cdot\|_r^{(P)} : 0<r<1\}$ of seminorms.
Clearly, for each $T\in\cB(H)^n$ with $\| T\|<1$
the functional calculus $\gamma_T^\free$ (see above) is continuous.

In fact, it easily follows from the results of \cite{Pop_holball} that
$\cF(\BB^n)$ admits the following simpler characterization.

\begin{prop}
\label{prop:freeball_series}
We have
\[
\cF(\BB^n)=
\Bigl\{ a=\sum_{\alpha\in W_n} c_\alpha\zeta_\alpha :
\| a\|_r=\sum_{k=0}^\infty\Bigl(\sum_{|\alpha|=k} |c_\alpha|^2\Bigr)^{1/2} r^k<\infty
\;\forall 0<r<1\Bigr\}.
\]
The families $\{\|\cdot\|_r : 0<r<1\}$ and $\{\|\cdot\|_r^{(P)} : 0<r<1\}$
of seminorms are equivalent.
\end{prop}

Proposition~\ref{prop:freeball_series} can also be interpreted as follows.
Let $\ell_1^n$ (respectively, $\ell_2^n$) denote the $n$-dimensional version
of $\ell_1$ (respectively, of $\ell_2$). Then Taylor's polydisk algebra $\cF_n(1)$
is a weighted $\ell_1$-sum of the projective tensor powers
$\ell_1^n\Tens_\pi\cdots\Tens_\pi\ell_1^n$, while $\cF(\BB^n)$ is
a weighted $\ell_1$-sum of the Hilbert tensor powers
$\ell_2^n\Tens_{\mathrm{hilb}}\cdots\Tens_{\mathrm{hilb}}\ell_2^n$.

Using Proposition~\ref{prop:freeball_series}, one can easily show that
$\cF(\BB^n)$ is a Schwartz space, but is not a nuclear space
(cf. Remark~\ref{rem:Tay_alg}). Thus $\cF(\BB^n)$ is not an HFG algebra.
For our purposes, $\cF(\BB^n)$ is useful because of the following fact.

\begin{prop}
\label{prop:q_ball_quot_free_ball}
\begin{mycompactenum}
\item For each $q\in (0,1)$, $\cO^\hol_q(\BB^n)$ is topologically isomorphic to the
quotient of $\cF(\BB^n)$ by the two-sided closed ideal generated by the elements
$\zeta_i\zeta_j-q\zeta_j\zeta_i$ for all $i<j$.
\item $\cO^\hol(\BB^n)$ is topologically isomorphic to the
quotient of $\cF(\BB^n)$ by the two-sided closed ideal generated by the elements
$\zeta_i\zeta_j-\zeta_j\zeta_i$ for all $i<j$.
\end{mycompactenum}
\end{prop}

Proposition~\ref{prop:q_ball_quot_free_ball} looks similar to Proposition~\ref{prop:q_poly_quot_free_poly},
but its proof is less elementary. The proof of part (ii) involves, in particular,
V.~M\"uller's characterizations~\cite{Muller_sprad} of the joint spectral radius
in commutative Banach algebras, as well as J.~L.~Taylor's holomorphic functional calculus on $\BB^n$.
The proof of part (i) is based on Theorem~\ref{thm:poly_vs_ball}.

Repeating now the construction of the bundle $(D,p)$ with $\cF(\DD_r^n)$ replaced by
$\cF(\BB^n)$, and using Proposition~\ref{prop:q_ball_quot_free_ball} instead of
Proposition~\ref{prop:q_poly_quot_free_poly}, we obtain the bundle $(B,p)$
whose existence was stated in part (ii) of Theorem~\ref{thm:bnd}.

\begin{ackn}
The author thanks A.~Ya.~Helemskii, D. Proskurin, and B. Solel for helpful discussions,
and the referee for a careful reading and useful comments.
\end{ackn}

\end{document}